\documentclass[12pt,reqno]{amsart}
\begin{document}
\numberwithin{equation}{section}


\def\1#1{\overline{#1}}
\def\2#1{\widetilde{#1}}
\def\3#1{\widehat{#1}}
\def\4#1{\mathbb{#1}}
\def\5#1{\frak{#1}}
\def\6#1{\mathcal{#1}}

\def\N{{\4N}}
\def\C{{\4C}}
\def\R{{\4R}}
\def\M{{\3M}}


\def\cn{{\C^n}}
\def\cnn{{\C^{n'}}}
\def\ocn{\2{\C^n}}
\def\ocnn{\2{\C^{n'}}}


\def\const{{\rm const}}
\def\rk{{\rm rank\,}}
\def\id{{\sf id}}
\def\Aut{{\sf Aut}}
\def\CR{{\rm CR}}
\def\Orb{{\sf Orb}}

\def\codim{{\rm codim}}
\def\crd{\dim_{{\rm CR}}}
\def\crc{{\rm codim_{CR}}}

\def\GL{{\bf GL}}
\def\phi{\varphi}
\def\eps{\varepsilon}
\def\d{\partial}
\def\a{\alpha}
\def\z{{\bar z}}
\def\l{\lambda}
\def\L{\Lambda}
\def\Im{\text{{\sf Im\,}}}
\def\Re{\text{{\sf Re\,}}}

\emergencystretch15pt \frenchspacing

\newtheorem{Thm}{Theorem}[section]
\newtheorem{Cor}[Thm]{Corollary}
\newtheorem{Prop}[Thm]{Proposition}
\newtheorem{Lem}[Thm]{Lemma}
\theoremstyle{definition}\newtheorem{Def}[Thm]{Definition}
\theoremstyle{remark}\newtheorem{Rem}[Thm]{Remark}
\newtheorem{Exa}[Thm]{Example}

\def\b{\begin}
\def\e{\end}

\def\bl{\b{Lem}}
\def\el{\e{Lem}}
\def\bp{\b{Prop}}
\def\ep{\e{Prop}}
\def\bt{\b{Thm}}
\def\et{\e{Thm}}
\def\bc{\b{Cor}}
\def\ec{\e{Cor}}
\def\bd{\b{Def}}
\def\ed{\e{Def}}
\def\br{\b{Rem}}
\def\er{\e{Rem}}
\def\bpr{\b{proof}}
\def\epr{\e{proof}}
\def\ben{\b{enumerate}}
\def\een{\e{enumerate}}
\def\be{\b{Eq}}
\def\ee{\end{equation}}

\def\Label#1{\label{#1}}

\def\bl{\begin{Lem}}
\def\el{\end{Lem}}
\def\bp{\begin{Prop}}
\def\ep{\end{Prop}}
\def\bt{\begin{Thm}}
\def\et{\end{Thm}}
\def\bc{\begin{Cor}}
\def\ec{\end{Cor}}
\def\bd{\begin{Def}}
\def\ed{\end{Def}}
\def\br{\begin{Rem}}
\def\er{\end{Rem}}
\def\be{\begin{Exa}}
\def\ee{\end{Exa}}
\def\bpf{\begin{proof}}
\def\epf{\end{proof}}
\def\ben{\begin{enumerate}}
\def\een{\end{enumerate}}

\title[Remarks on the rigidity of CR-manifolds]{Remarks on the rigidity of CR-manifolds}
\author[S. Kim and D. Zaitsev]{Sung-Yeon Kim and Dmitri Zaitsev}
\address{Korea Institute for Advanced Study, Cheongnyangni 2-dong Dongdaemun-gu Seoul
130-722, Korea } \email{sykim87@kias.re.kr}
\address{School of Mathematics, Trinity College Dublin, Dublin 2, Ireland}
\email{zaitsev@maths.tcd.ie} 

\begin{abstract}
We propose a procedure to construct new smooth CR-manifolds
whose local stability groups, equipped with their natural topologies, 
are subgroups of certain (finite-dimensional) Lie groups but not Lie groups themselves.
\end{abstract}

\maketitle

\section{Introduction}

Given a germ $(M,p)$ of a real submanifold of $\C^n$, its basic
invariant is the {\em local stability group} $\Aut(M,p)$, i.e.\ the group of
all germs at $p$ of local biholomorphic maps of $\C^n$ fixing $p$ and preserving
the germ $(M,p)$.
By the work of several authors \cite{CM,BERasian,Z97,ELZ,LM}
it is known that this group is a (finite-dimensional) Lie group
(in the natural inductive limit topology)
for germs of {\em real-analytic} submanifolds satisfying certain nondegeneracy conditions,
e.g.\ those having nondegenerate Levi form.
On the other hand, in the absence of the nondegeneracy conditions,
the group $\Aut(M,p)$ can possibly be infinite-dimensional
(in the sense that it contains Lie groups of arbitrarily large dimension).
(E.g.\ the local stability group of $(\R,0)$ in $\C$ consists
of all convergent power series with real coefficients.)
Furthermore, recent results in \cite{BRWZ} show that a similar principle also holds
for {\em global CR-automorphisms}, both real-analytic and smooth.

One purpose of this paper is to show that, in contrast with the behaviour mentioned above, 
a similar alternative does not anymore hold for the {\em local 
stability group} of a {\em smooth} real submanifold.
In particular, we show that, for any $n\ge 2$, there
exists a germ $(M,p)$ of a smooth strongly pseudoconvex hypersurface
in $\C^n$ with $\Aut(M,p)$ being (topologically) isomorphic to a countable
dense subgroup of the circle $S^1\subset\C$ and hence not being a Lie group. 
In fact, $\Aut(M,p)$ can be arranged to be isomorphic to
any increasing countable union of finite subgroups of $S^1$,
for instance, to the subgroup
\begin{equation}\Label{dense-sub}
\{e^{2\pi i \frac{ l}{2^m}} : l,m\in\N \} \subset S^1.
\end{equation}

Furthermore, our construction yields similar properties also for the (generally larger) 
{\em local CR stability group} $\Aut_\CR(M,p)$,
consisting of all germs at a point $p$ of smooth CR-automorphisms of $M$ fixing $p$.
Recall that a germ of a smooth transformation $\phi\colon (M,p)\to (M,p)$
is a {\em CR-automorphism} if it preserves the subbundle $T^cM$ and
and the restriction of its differential $d\phi|_{T^cM}$ is $\C$-linear, where
$$T^cM := TM\cap iTM.$$

Another purpose of this paper is to provide a general construction
of new smooth generic submanifolds with certain prescribed local CR stability groups
(recall that a real submanifold $M\subset\C^n$ is {\em generic} if $T_qM + iT_qM =T_q\C^n$
for all $q\in M$).
More precisely, we show the following:

\begin{samepage}
\bt\Label{group}
Let $(M,p)$ be a germ of a smooth generic submanifold in $\cn$ of
positive codimension and of finite type and assume
that it is invariant under an increasing countable union $G$ of finite subgroups of $\Aut(\cn,p)$.
Then there exists a $G$-invariant germ of another smooth generic
submanifold $(\2M,p)$ of the same dimension as $(M,p)$, which is tangent to $(M,p)$
of infinite order and has the following properties:
\begin{enumerate}
\item[(i)] $\Aut(\2M,p)=G$;
\item[(ii)] $\Aut_\CR(\2M,p)=\{g|_M : g\in G\}$.
\end{enumerate}
\et
\end{samepage}

We use here the notion of  finite type due to {\sc Kohn}~\cite{K72} and
{\sc Bloom-Graham}~\cite{BG}): a germ $(M,p)$ is of {\em finite type},
if all germs at $p$ of smooth real vector fields on $M$ tangent to $T^cM$
span together with their iterated commutators the full tangent space $T_pM$.

We now illustrate Theorem~\ref{group} by an example, where it can be applied.

\be
Consider a real hypersurface $M\subset\C^{n+1}$, $n\ge 1$, 
given in coordinates $(z_1,\ldots,z_n,w)\in\C^{n+1}$ by
$$\Im w =\phi (|z_1|^2,z_2,\ldots,z_n,\bar z_2,\ldots,\bar z_n,\Re w),$$
where $\phi$ is any smooth function such that $0\in M$ and $(M,0)$ is of finite type.
Then $(M,0)$ is clearly invariant under the rotation group consisting of all transformations
$(z_1,z_2,\ldots,z_n,w)\mapsto (e^{2\pi i\theta} z_1,z_2,\ldots,z_n,w)$ for all real $\theta$.
Now we can take the subgroup $G$ consisting of all these transformations corresponding
to $\theta = l/2^m$ with $l,m$ being positive integers.
Then $G$ clearly satisfies the assumptions of Theorem~\ref{group}.
We then conclude that there exists a new real submanifold $\2M\subset\C^{n+1}$
such that both $\Aut(\2M,p)$ and $\Aut_\CR(\2M,p)$ are (topologically) isomorphic to $G$,
which is a topological subgroup of $S^1$ but is not itself a Lie group.
Similar examples can be obtained for other $S^1$-actions
or $S^1\times S^1$-actions or actions by more general compact groups
leaving $(M,p)$ invariant, where $M$ can also be of any codimension.
\ee

As a remarkable consequence of Theorem~\ref{group}
and the mentioned results \cite{BERasian,Z97}, our construction 
provides germs of smooth generic submanifolds (even of
strongly pseudoconvex hypersurfaces) that are not CR-equivalent to
any germ of any real-analytic CR-manifold. 
Recall that $(M,p)$ is called {\em finitely nondegenerate} if
\begin{equation}\Label{span}
{\sf span}_\C\left\{
L_1\ldots L_k \rho^j_Z(p) : k\ge0, \, 1\le j\le d
\right\}=\cn,
\end{equation}
where $\rho=(\rho^1,\ldots,\rho^d)$ is a defining function of $M$ near $p$
with $\partial\rho^1\wedge\cdots\wedge\partial\rho^d\ne0$,
$\rho^j_Z$ denotes the complex gradient of $\rho^j$ in $\cn$
and the span is taken over all collections of
germs at $p$ of smooth $(0,1)$ vector fields $L_1\ldots L_k$ on $M$.
We have:

\bc For any germ $(M,p)$ of a smooth generic submanifold in $\cn$,
which is finitely nondegenerate and of finite type,
and which is invariant under the group $S^1:=\{z\in\C : |z|=1\}$ acting on
$\cn$ by multiplication, there exists a germ of another smooth
generic submanifold $(\2M,p)$, tangent to $(M,p)$ of infinite
order, which is not CR-equivalent to any germ of a real-analytic
generic submanifold of $\cn$. \ec

Indeed, since $S^1$ has many dense subgroups $G$ satisfying the assumptions of Theorem~\ref{group}
(e.g.\ the subgroup in \eqref{dense-sub}),
Theorem~\ref{group} yields a germ $(\2M,p)$ whose local CR stability group
is isomorphic to $G$ and hence is not topologically isomorphic to any Lie group
(and not locally compact).
On the other hand, the local CR stability group
of any real-analytic generic submanifold of $\cn$, which is CR-equivalent to $(\2M,p)$
(and hence is also finitely nondegenerate and of finite type),
is known to be always a Lie group (see \cite{BERasian} for hypersurfaces and \cite{Z97} 
for higher codimension). Since the local CR stability group is a CR-invariant,
$(\2M,p)$ cannot be CR-equivalent to any germ of a real-analytic
generic submanifold of $\cn$.

{\bf Acknowledgement.}
The authors would like to thank the referee for careful reading of
the manuscript and many critical suggestions.

\section{Jet spaces and jet groups}
Here we recall the jet terminology and introduce the notation 
that will be used throughout the paper.
Recall that, given two complex manifolds $X$ and $X'$ and an
integer $k\ge 0$, a {\em $k$-jet of a holomorphic map} is an equivalence
class of holomorphic maps from open neighborhoods of $x$ in $X$ into
$X'$ with fixed partial derivatives at $x$ up to order $k$. Denote
by $J^k_x(X,X')$ the set of all such $k$-jets. The union
$J^k(X,X'):=\bigcup_{x\in X}J^k_x(X,X')$ carries a natural fiber
bundle structure over $X$. For a holomorphic map $f$ from a
neighborhood of $x$ in $X$ into $X'$, denote by $j^k_xf\in
J^k_x(X,X')$ the corresponding $k$-jet. 
In local coordinates, $j^k_xf$ can be represented by
the coordinates of the reference point $x$
and all partial derivatives of $f$ at $x$ up to order $k$.
If $X$ and $X'$ are smooth
algebraic varieties, $J^k_x(X,X')$ and $J^k(X,X')$ are also of this type. 
We also denote by $J^k_{x,x'}(X,X')$ the space of all
$k$-jets sending $x$ into $x'$. The subset $G^k_x(X)\subset
J^k_{x,x}(X,X)$ of all {\em invertible} $k$-jets forms an
algebraic group with respect to composition.
Completely analogously {\em $k$-jets of smooth maps} between smooth real manifolds
$M$ and $M'$ are defined, for which we shall use the same notation $J^k(M,M')$.
The possible confusion will be eliminated by
the convention that we write $X_\R$
whenever we consider a complex manifold $X$ as a real manifold.
Thus, if $X$ and $X'$ are complex manifolds,
$J^k(X,X')$ is the space of all $k$-jets of holomorphic maps
and $J^k(X_\R,X'_\R)$ is the space of all $k$-jets of smooth maps.

Furthermore, we shall need $k$-jets of real submanifolds of fixed dimension of
a smooth real manifold $M$. Let
$\6C_x^{k,m}(M)$ be the set of all germs at $x$ of real
$C^k$-smooth $m$-dimensional submanifolds of $M$ through $x$. We
say that two germs $V,V'\in\6C_x^{k,m}(M)$ are {\em $k$-equivalent},
if, in a local coordinate neighborhood of $x$ of the form
$U_1\times U_2$, $V$ and $V'$ can be given as graphs of smooth
maps $\phi,\phi'\colon U_1\to U_2$ such that
$j^k_{x_1}\phi=j^k_{x_1}\phi'$, where $x=(x_1,x_2)\in U_1\times
U_2$. Denote by $J_x^{k,m}(M)$ the set of all $k$-equivalence
classes of $\6C_x^{k,m}(M)$ and by $J^{k,m}(M)$ the union
 $\bigcup_{x\in M}J^{k,m}_x(M)$ with the natural fiber bundle structure over
$M$. Furthermore, for any real $C^k$-smooth $m$-dimensional
submanifold $V\subset M$ through $x$, denote by $j^k_x(V)\in
J_x^{k,m}(M)$ the corresponding $k$-jet. The space
$J_x^{k,m}(M)$ carries a natural real (nonsingular) algebraic variety
structure. 

We now introduce the notions of equivalence and rigidity 
that will be crucial in the sequel.
\bd\Label{rig}
\begin{enumerate}
\item Two $k$-jets of real submanifolds of the same dimension 
$\L_j\in J^{k,m}_{p_j}(\C^n_\R)$, $j=1,2$, are called
{\em biholomorphically equivalent} if there exists a germ of a
biholomorphic map $(\cn,p_1)\to (\cn,p_2)$ sending $\L_1$ to
$\L_2$.
\item A $C^k$-smooth generic submanifold $M\subset\cn$ is called 
{\em totally rigid of order $k$},
if for any $p_1\ne p_2\in M$, the jets $j^k_{p_1}(M)$ and
$j^k_{p_2}(M)$ are not biholomorphically equivalent in the sense of (1).
\item A $k$-jet $\L\in J^{k,m}_p(\C^n_\R)$ is called {\em totally rigid}
if any $C^k$-smooth submanifold passing through $p$ and having
$\L$ as its $k$-jet at $p$, contains a neighborhood of $p$ that is
totally rigid of order $k-1$ in the sense of (2).
\end{enumerate}
\ed

\be
Any $0$-jets of real submanifolds at $p$ are obviously biholomorphically equivalent
and $1$-jets are equivalent if and only if their CR-dimensions at $p$ are the same.
Two $2$-jets of generic submanifolds are equivalent if and only their Levi forms at $p$
are linearly equivalent (e.g.\ of the same rank and signature in the hypersurface case).
Furthermore it follows from the Chern-Moser theory \cite{CM} that two
$k$-jets of Levi-nondegenerate hypersurfaces of the same signature
are always biholomorphically equivalent for $k\le 5$ in case $n=2$ and for $k\le 3$ in case $n>2$,
but may not be equivalent in general for $k$ larger.
\ee

It is crucial for our method to consider the total rigidity of order $k-1$ in (3)
(rather then e.g.\ of order $k$) for the representing submanifolds $M$ with $j^k_{p}(M)=\L$.
This allows us to achieve the total rigidity of $M$ of order $k-1$ near $p$
by ensuring that the first order derivatives of $j^{k-1}_q(M)$ at $p$ as function of $q\in M$
have suitable transversality property with respect to the submanifolds (orbits)
of biholomorphically equivalent $(k-1)$-jets
(see the proof of Proposition~\ref{rigid} below for more details).
Thus we need $\L$ to be of higher order than $k-1$ to include the extra derivatives.
More precisely, the existence of totally rigid jets is guaranteed by the following
statement.

\bp\Label{rigid} For fixed integers $n<m<2n$, a point $p\in\cn$ and
sufficiently large $k$ (depending on $n$ and $m$ but not on $p$), 
the set of all totally rigid $k$-jets in
$J^{k,m}_p(\C^n_\R)$ contains an open dense subset. 
\ep

In fact we show that  the number $k$ in Proposition~\ref{rigid} can be
chosen such that the following inequality holds:
\begin{equation}\Label{ineq}
(2n-m)\bigg(\Big(\begin{matrix} k+m-1 \\ m \end{matrix}\Big)
-1\bigg) - 2n\bigg(\Big(\begin{matrix} k+n-1 \\ n
\end{matrix}\Big) -1\bigg) \ge m.
\end{equation}

The proof will be based on the following lemmas
(of which the first is standard and provided with a short proof for the reader's convenience).
We write $\|\cdot\|_{C^l}$ for the standard $C^l$ norm.

\bl\Label{perturb} Let $\phi\colon\R^n\to\R^m$ be a smooth map and
$V$ be an open neighborhood in $\R^n$ of a point $a\in \R^n$. Then
for any $\eps>0$ and any integers $0\le k\le l$, there exists $\delta>0$
such that, if $\L\in J^k_a(\R^n,\R^m)$ is a $k$-jet with
$\|\L-j^k_a \phi\|<\delta$, then there exists another smooth map
$\2\phi\colon\R^n\to \R^m$ such that $j^k_a\2\phi=\L$,
$\2\phi|_{\R^n\setminus V} =  \phi|_{\R^n\setminus V}$ and
$\|\2\phi-\phi\|_{C^l}<\eps$.
 \el

\bpf Without loss of generality, $V$ is bounded. We shall look for
a map $\2\phi\colon\R^n\to \R^m$ of the form
\begin{equation}\Label{psi}
\2\phi(x):=\phi(x)+\chi(x)\cdot (\psi(x)-\phi(x)),
\end{equation}
where $\psi\colon \R^n\to\R^m$ is a smooth map with $j^k_a\psi=\L$
and $\chi\colon\R^n\to\R$ is a fixed smooth function which is $1$
in a neighborhood of $a$ and $0$ outside $V$. Then
$j^k_a\2\phi=\L$ and $\2\phi|_{\R^n\setminus V} =
\phi|_{\R^n\setminus V}$. Furthermore, there exists $C>0$
depending on $\chi$ but not on $\psi$ such that
$$\|\2\phi-\phi\|_{C^l}<C \|\psi-\phi\|_{C^l}.$$
If $\delta$ is
sufficiently small and $\L$ satisfies our assumption, we can always choose $\psi$
with $\|\psi-\phi\|_{C^l} < \eps/C$ on the closure $\1V$. Then the
map $\2\phi$ given by \eqref{psi} satisfies the required
properties. \epf

\bl\Label{transversal}
Let $k,l,s\ge0$ and $0\le r\le m$ be any integers, $U\subset J^l(\R^m,\R^s)$ be an open set
and $F\colon U\to \R^r$ be a smooth map of constant rank $r$.
Then, for any nonempty open set $B\subset \R^m$ and a smooth map $f\colon B\to \R^s$ 
satisfying $j^l_xf\in U$ for all $x\in B$, there exists another nonempty open subset $\3B\subset B$
and another smooth map $\3f\colon \3B\to \R^s$ such that the following holds:
\begin{enumerate}
\item the $C^k$ norm of $\3f-(f|_{\3B})$ can be chosen arbitrarily small;
\item $j^l_x\3f\in U$ for all $x\in \3B$;
\item the map $x\in \3B \mapsto F(j^l_x\3f)\in \R^r$ is also of constant rank $r$.
\end{enumerate}
\el

\bpf
Without loss of generality, the integer $k$ (used only in (1)) is $\ge l$.
We prove the lemma by induction on $r$. For $r=0$ the statement is trivial.
Suppose it holds for any $r<r_0\le m$ and we are given a map $F\colon U\to \R^{r_0}$
as in the lemma. Consider the standard splitting $\R^{r_0}=\R^{r_0-1}\times\R$
and the corresponding components $F_1\colon U\to\R^{r_0-1}$
and $F_2\colon U\to \R$ of $F$ (so that $F=(F_1,F_2)$). Then the induction assumption for $F_1$
yields a map $\3f\colon \3B\to \R^s$ such that
$\|\3f-(f|_{\3B})\|_{C^k}$ is arbitrarily small,
 $j^l_x\3f\in U$ for all $x\in \3B$ and the map
$\Phi_1 (x):= F_1(j^l_x\3f)\in \R^{r_0-1}$ is of constant rank $r_0-1<m$ for $x\in \3B$.
By shrinking $\3B$ if necessary, we may assume it is connected and of the form
$\3B=\3B_1\times \3B_2\subset \R^{m-r_0+1}\times\R^{r_0-1}$
and such that, for some $c\in \R^{r_0-1}$, the level set 
$C:=\{x\in \3B : \Phi_1 (x)=c\}$ is a graph of a smooth map $\phi\colon \3B_1\to \3B_2$.

We now consider two points $x_1\ne x_2\in C$ and look for a small perturbation $\2f$ of $\3f$
and $\2x_2$ of $x_2$ such that $x_1$ and $\2x_2$ still belong to the same level set $\2C$ of 
$\2\Phi_1 (x):= F_1(j^l_x\2f)$ but $F_2(j^l_{x_1}\2f) \ne F_2(j^l_{\2x_2}\2f)$.
More precisely, suppose that $F(j^l_{x_1}\3f) = F(j^l_{x_2}\3f)$
(otherwise no perturbation is needed).
By the assumption of the lemma, the map $F=(F_1,F_2)$ has constant rank $r_0$ in $U$.
Since $r_0\le m\le \dim U$, we can find a point $\2x_2\in \3B$
and a jet $\L\in U\cap J^l_{\2x_2}(\R^m,\R^s)$ arbitrarily close to $j^l_{x_2}\3f$
such that $F_1(\L)=F_1(j^l_{x_2}\3f)$ but $F_2(\L)\ne F_2(j^l_{x_2}\3f)$.
Since $\L$ can be chosen arbitrarily close to $j^l_{x_2}\3f$,
it can be represented by a smooth map $\2f$ that is arbitrarily close to $\3f$
in the $C^k$ norm and differs from it on a neighborhood of $x_2$ with compact support in 
$\3B\setminus\{x_1\}$ in view of Lemma~\ref{perturb}.
By choosing the norm $\|\2f-\3f\|_{C^k}$ sufficiently small,
we shall preserve the properties that $j^l_x\2f\in U$ for all $x\in \3B$,
the rank of $\2\Phi_1(x):= F_1(j^l_x\2f)$ is still $r_0-1$ and
 the level set $\2C:=\{x\in \3B : \2\Phi_1 (x)=c\}$
is still a graph of a smooth map $\2\phi\colon \3B_1\to \3B_2$.
Hence $\2C$ is a connected manifold containing two points $x_1,\2x_2\in \2C$
such that $F_2(j^l_{x_1}\2f)\ne F_2(j^l_{\2x_2}\2f)$.
The latter fact implies that the function $\2\Phi_2(x):=F_2(j^l_x\2f)$
is not constant on $\2C$ 
and therefore its differential is somewhere nonzero.
Putting this property together with the rank property of $\2\Phi_1$,
we conclude that the rank of $\2\Phi(x):=F(j^l_x\2f)$ is $r_0$ at some point
$x_0\in \3B$. The required conclusion is obtained by replacing $\3f$ with $\2f$
and $\3B$ with a sufficiently small open neighborhood of $x_0$.
\epf

\bpf[Proof of Proposition~\ref{rigid}]
We may assume $p=0$.
Consider the natural action of the group $G^{k-1}_0(\cn)$ 
(consisting of all $(k-1)$-jets at $0$ of local biholomorphic maps of $\cn$)
on the space $J^{k-1,m}_0(\C^n_\R)$ 
(consisting of all $(k-1)$-jets at $0$ of real $m$-dimensional
submanifolds of $\C^n_\R$ passing through $0$. 
The dimensions of the jet spaces
can be computed directly:
\begin{equation}\Label{gj}
\dim_\R G^{k-1}_0(\cn)= 2n\bigg(\Big(\begin{matrix} k+n-1 \\ n
\end{matrix}\Big) -1\bigg), \quad \dim_\R J^{k-1,m}_0(\C^n_\R) =
(2n-m)\bigg(\Big(\begin{matrix} k+m-1 \\ m \end{matrix}\Big)
-1\bigg).
\end{equation}
Hence the inequality (\ref{ineq}) is equivalent to
\begin{equation}\Label{codim}
\dim_\R J^{k-1,m}_0(\C^n_\R) - \dim_\R G^{k-1}_0(\cn) \ge m.
\end{equation}
In particular, for any sufficiently large $k$,
all orbits of $G^{k-1}_0(\cn)$ in $J^{k-1,m}_0(\C^n_\R)$ have their (real)
codimension at least $m$. In the rest of the proof we shall assume that
(\ref{codim}) is satisfied.

It is easy to see that $G^{k-1}_0(\cn)$ is an algebraic group acting rationally on $J^{k-1,m}_0(\C^n_\R)$
by calculating the group operation and the action in local coordinates.
Hence the orbits of $G^{k-1}_0(\cn)$ form,
on an open dense subset $\Omega\subset J^{k-1,m}_0(\C^n_\R)$, 
a foliation into real submanifolds of a fixed constant codimension $\ge m$.
Consider any $k$-jet $\L_0\in J^{k,m}_0(\C^n_\R)$ represented by
the graph of a $C^\infty$-smooth map
$\phi_0\colon\R^m\to\R^{2n-m}$ with $\phi_0(0)=0$,
where we choose a suitable identification of $\C^n_\R$ with $\R^m\times\R^{2n-m}$
(after a possible permutation of the real coordinates).
Thus $\L_0=j^k_0 ((\id\times\phi_0)(\R^m))$.
By the density of $\Omega$, we can find
another $C^\infty$-smooth map $\phi\colon\R^m\to\R^{2n-m}$
with $\phi(0)=0$ and $\Theta:=j^k_0\phi$ arbitrarily close to $j^k_0\phi_0$
such that the $(k-1)$-jet $\L\in J^{k-1,m}_0(\C^n_\R)$ 
at $0$ of the graph of $\phi$ is contained in $\Omega$.
(By this choice, also $\L$ is arbitrarily close to $\L_0$.)
We can now find an open neighborhood $U$ of $j^{k-1}_0\phi$ in $J^{k-1}(\R^m,\R^{2n-m})$ 
and a smooth map $F\colon U\to \R^m$ of constant rank $m$ and constant on the orbits
(recall that $m$ does not exceed the orbit codimension),
such that two $(k-1)$-jets $\L_j\in J^{k-1,m}_{(x_j,x'_j)}(\C^n_\R)$, $j=1,2$, near $\L_0$
(where $(x_j,x'_j)\in\R^m\times\R^{2n-m}$),
which are represented by graphs of some smooth maps $\phi_1,\phi_2\colon \R^m\to\R^{2n-m}$, 
are biholomorphically inequivalent (in the sense of Definition~\ref{rig})
whenever $F(j^{k-1}_{x_1}\phi_1) \ne F(j^{k-1}_{x_2}\phi_2)$.
Here $F$ can be obtained by taking the first $m$ coordinates
in any real coordinate system $(x,y)\in \R^{m_1}\times \R^{m_2}$,
for which the orbits are given by $x=\const$.
Then Lemma~\ref{transversal} can be applied to $U$, $F$, $f:=\phi$ and 
an arbitrarily small neighborhood of $B\subset \R^m$ of $0$.
Let $\3B$ and $\3f\colon \3B\to \R^{2n-m}$ be given by the lemma.

We claim that, for any $x_0\in\3B$, the $k$-jet $\L(x_0)\in J^{k,m}(\C^n_\R)$
of the graph of $\3f$ at $(x_0,\3f(x_0))$ is totally rigid in the sense of Definition~\ref{rig} (3).
Indeed, fix any $x_0\in\3B$ and consider any $C^k$-smooth real 
$m$-dimensional submanifold $V\subset\C^n_\R$ passing through $(x_0,\3f(x_0))$
with $j^k_{(x_0,\3f(x_0))}(V) = \L(x_0)$.
By shrinking $V$, if necessary, we may assume that $V$ is a graph of 
a smooth map $g\colon B(x_0)\to \R^{2n-m}$, where $B(x_0)$ is a suitable open neighborhood of $x_0$.
Then $j^k_{x_0} g=j^k_{x_0}\3f$ and therefore, 
the ranks of the maps $x\mapsto F(j^{k-1}_x g)$ and $x\mapsto F(j^{k-1}_x \3f)$
coincide at $x=x_0$. (Here is the step, where we use the different integers $k$ and $k-1$
for $x_0$ and points nearby respectively.)
By property (3) in Lemma~\ref{transversal},
the rank of the second map is $m$ and hence, so is the rank of the first map.
But the latter fact implies that $F(j^{k-1}_{x_1} g)\ne F(j^{k-1}_{x_2} g)$
for any $x_1\ne x_2$ sufficiently close to $x_0$.
In view of the choice of $F$, it follows that the jets 
$j^{k-1}_{(x_1,g(x_1))}(V)$ and $j^{k-1}_{(x_2,g(x_2))}(V)$
are biholomorphically inequivalent for any such $x_1\ne x_2$,
which is precisely what is needed to show that $\L(x_0)$ is totally rigid.
It remains to observe, that any translation of $\L(x_0)$ is also totally rigid,
hence we can find totally rigid $k$-jets also in $J^{k,m}_0(\C^n_\R)$  arbitrarily close to the original 
jet $\L_0$. The proof is complete. \epf

\section{Realization of certain groups as CR stability groups}

In the sequel we shall use the same letter for a germ and its
representative unless there will be a danger of confusion. We
begin with a standard lemma, whose proof is given here for the
reader's convenience.

\bl\Label{sequence} Let $(G_k)_{k\ge 1}$ be an increasing sequence
of finite groups of germs of local biholomorphic maps of  $\cn$
in a neighborhood of a point $p$ fixing that point.
Let $M$ be a smooth real submanifold of $\cn$ passing through $p$.
Let $(D_k)_{k\ge 1}$ be a sequence of domains containing $p$
and such that, for each $k$, the germs from $G_k$ can be represented
by biholomorphic self-maps of $D_k$. Then there exist a sequence of points
$p_k\in M$, $k\ge 1$, converging to $p$ and a sequence of mutually disjoint open
neighborhoods $V_k$ of $p_k$ in $\cn$ such that $\1V_k\subset D_k$
and, if $g(\overline V_k)\cap \overline V_l\ne\emptyset$ for some
$k$, $l$ and $g\in G_k$, then necessarily $k=l$ and $g\equiv\id$. \el

Note that the existence of domains $D_k$ easily follows from the
finiteness of each $G_k$. Indeed, if $\2D_k$ is any domain
where all germs from $G_k$  biholomorphically extend,
it suffices to take $D_k:= \bigcap_{g\in G_k} g(\2D_k)$.

\bpf We shall construct $p_k\in M$ and $V_k\subset\cn$ inductively.
Let $k=1$. Since $G_1$ is finite, the set of points $x\in D_1$,
such that there are two elements $g_1\ne g_2\in G_1$ with $g_1(x)=
g_2(x)$, is a complement of a proper analytic subset. Hence we can
choose $p_1\in M$ and a neighborhood $V_1\subset\subset D_1$ of
$p_1$ in $\cn$ with $p\not\in\1V_1$ such that  $g(\overline
V_1)\cap \overline V_1\ne \emptyset$ for $g\in G_1$ implies
$g\equiv\id$.

Now suppose that $p_l$ and $V_l$ with $p\not\in\1V_l$ have been chosen for
all $l< k$. Since $G_k$ is finite, we can choose a neighborhood
$U$ of $p$ in $D_k$ such that $g(U)\cap \overline V_l=\emptyset$
for all $l<k$ and $g\in G_k$. Using the same argument as before we
can choose $p_k$ arbitrarily close to $p$ and $V_k$ with $p\notin\1V_k$ such that
$p_k\in V_k \subset \1V_k\subset U$ and, for any $g\not\equiv\id\in G_k$, $g(\overline
V_k)\cap \overline V_k=\emptyset$. Since $\1V_k\subset U$ and
$g(U)\cap \overline V_l=\emptyset$ for all $l<k$ and $g\in G_k$,
it follows that $g(\overline V_k)\cap \overline V_l\neq\emptyset$ can hold for some
$l\le k$ and $g\in G_k (\supset G_l)$ if and only if $k=l$ and $g\equiv\id$. It is
easy to see that the sequences $(p_k)$ and $(V_k)$ so constructed
satisfy the required properties. \epf

\bpf[Proof of Theorem~{\rm\ref{group}}]
Since the statement is local, we fix 
an identification $\C^n_\R\cong \R^m\times\R^{2n-m}$ near $p$
such that $M$ is represented by the graph of a smooth map $\phi\colon\R^m\to\R^{2n-m}$
with $\|\phi\|_{C^1}$ sufficiently small.
We shall write $B_r(a)$ for the open ball with center $a$ and radius $r$
with respect to the product metric of the Euclidean metrics of $\R^m$ and $\R^{2n-m}$.
With this choice of the metric on $\R^m\times\R^{2n-m}$ and $\|\phi\|_{C^1}$ sufficiently small,
we have the property that, for any $a\in M$, the intersection 
$B_r(a)\cap M$ coincides with the graph of $\phi$
over the projection of $B_r(a)$ to $\R^m$ (which is an Euclidean ball in $\R^m$).
In the course of the proof we shall consider small perturbations of $M$
obtained as graphs of small perturbations of $\phi$.
We shall always assume that the $C^1$ norms of these perturbations are still small,
so that the mentioned relation between
ball intersections with their graphs and graphs over balls still holds.

By the assumption, $G=\bigcup_k G_k$, where $G_k$, $k\ge 1$, is an
increasing sequence of finite groups of local biholomorphic maps of $\cn$
in a neighborhood of $p$, fixing $p$ and preserving the germ $(M,p)$.
As indicated above, we can choose a decreasing sequence of open neighborhoods $D_k$ of $p$ in
the unit ball $B_1(p)$ in $\cn$ centered at $p$, such that, for each $k$,
all germs in $G_k$ can be represented by biholomorphic self-maps of a neighborhood
of $\1D_k$. In the sequel we shall identify the germs from $G_k$ with
their biholomorphic representatives.

Let $p_k\in D_k\cap M$ and $V_k$ be given by Lemma~\ref{sequence}.
Moreover, we can make the above choice of $p_k$, $D_k$, $V_k$
inductively such that, in addition,
\begin{equation}\Label{supinf}
\max\Big(\sup_{z\in D_{k+1}, g\in G_{k+1}}\|g'(z)\| , \sup_{z\in D_{k}, g\in G_{k}}\|g'(z)\| \Big)
\frac{ \sup_{z\in V_{k+1}} |z-p| } {\inf_{z\in V_{k}} |z-p|} \to 0, \quad k\to\infty,
\end{equation}
where $g'$ is the Jacobian matrix of $g$.

For an $l$-jet $\L\in J^{l,m}_x(\C^n_\R)$, we denote by
$\Orb(\L)\subset J^{l,m}(\C^n_\R)$ the set of all $l$-jets
$\2\L\in J^{l,m}_y(\C^n_\R)$ for all $y\in\C^n_\R$ that are
biholomorphically equivalent to $\L$ (in the sense of
Definition~\ref{rig}). Then it follows from \eqref{gj}
that, for $l$ sufficiently large, the subset 
$\bigcup_{x\in M\cap D_1}\Orb(j^l_xM)\subset J^{l,m}(\C^n_\R)$ 
has Lebesgue measure zero and therefore its complement is dense in $J^{l,m}(\C^n_\R)$.
The same argument obviously applies to any other real submanifold of $\cn$
of the same dimension as $M$.
We shall use this property to choose jets that are not
in the certain unions of orbits.
We shall consider $l$ sufficiently large so that this choice is always possible.

We next consider a sequence $\eps_k$, $0<\eps_k<1$, converging to $0$ and such that
\begin{equation}\Label{eta}
\eps_k \Big(\sup_{x\in D_k, g\in G_k}\|g'(x)\| \Big)
\to 0, \quad k\to\infty.
\end{equation}
Then, as a consequence of
Lemma~\ref{perturb} and Proposition~\ref{rigid}, 
we can find a sequence of neighborhoods
$U_k$ of $p_k$ in $V_k$ with $\1U_k\subset V_k$ and a sequence of
graphs $N_k$ of smooth maps $\phi_k\colon \R^m\to\R^{2n-m}$ with
\begin{equation}\Label{epsk}
\|\phi_k - \phi\|_{C^k} < \eps_k, \quad N_k\setminus \1U_k=M\setminus \1U_k,
\end{equation}
and such that the following holds.
There exist points $q_k\in N_k$
and real numbers $\delta_k>0$ such that, for each $k$, the
$2\delta_k$-neighborhood of $q_k$ in $N_k$, 
$B_{2\delta_k}(q_k)\cap N_k$, is contained in $U_k$,
totally rigid (in  the sense of Definition~\ref{rig}) and
\begin{equation}\Label{orb-cond}
j^l_{q_k} N_k \notin \bigcup_{x\in M\cap D_1}\Orb(j^l_xM).
\end{equation}

As the next step, we define $X_k\subset N_k \cap U_k$ to be the subset of all
points $y\ne q_k$ for which there exists a CR-diffeomorphism
from $B_{\eps_k\delta_k}(y)\cap N_k$ into $B_{\delta_k}(q_k)\cap N_k$,
sending $y$ into $q_k$. 
Since the CR-manifold $B_{2\delta_k}(q_k)\cap N_k$ is totally rigid
by our construction, it is clear that for any $y_1\ne y_2\in X_k$,
the neighborhood $B_{\eps_k\delta_k}(y_1)\cap N_k$
cannot contain $y_2$.
Hence the neighborhoods $B_{\frac{\eps_k\delta_k}2}(y)\cap N_k$, $y\in X_k$, do not intersect
and so $X_k$ must be a finite set.
It is also clear that $X_k\cap B_{2\delta_k}(q_k)\cap N_k = \emptyset$
again by the total rigidity of $B_{2\delta_k}(q_k)\cap N_k$.
Furthermore we must have $X_k\subset N_k\cap U_k$ in view of \eqref{epsk} and \eqref{orb-cond}. 

We next choose a sequence $\eta_k$, $0<\eta_k<{\eps_k\delta_k}$
and apply again Lemma~\ref{perturb} to
obtain a sequence of graphs $M_k$ of smooth maps 
$\psi_k\colon \R^m\to \R^{2n-m}$ with
\begin{equation}\Label{epsk1}
\|\psi_k - \phi_k\|_{C^k} < \eps_k
\end{equation} 
and finite subsets $\2X_k\subset M_k$ with
\begin{equation}\Label{orb-cond1}
j^l_{y} M_k \notin \bigcup_{x\in N_k\cap D_1}\Orb(j^l_x N_k) \quad \forall y\in \2X_k,
\end{equation}
and such that, for $W_k:=\bigcup_{y\in \2X_k} B_{\eta_k}(y)$,
\begin{equation}\Label{mk-nk}
X_k\cap M_k = \emptyset, \quad
M_k\setminus \1{W_k} = N_k \setminus \1{W_k}, \quad \1{W_k}\cap B_{\delta_k}(q_k) = \emptyset.
\end{equation}
We may in addition assume
that $\eta_k$ is sufficiently small so that $\1W_k\subset V_k$.

Finally, we define the new generic submanifold $\2M\subset\C^n$ by replacing $g(M\cap
V_k)$ with $g(M_k\cap V_k)$ for every sufficiently large $k$ and
every $g\in G_k$, i.e.\
\begin{equation}\Label{consm}
\2M:=\Big(M\setminus \bigcup_{k\ge k_0, g\in G_k} g(M\cap
V_k)\Big)\bigcup \Big(\bigcup_{k\ge k_0,g\in G_k} g(M_k \cap
V_k)\Big),
\end{equation}
where $k_0$ is sufficiently large. 
(Note that all neighborhoods $g(V_k)$, $g\in G_k$, $k\ge k_0$,
are disjoint together with their closures since they are given by Lemma~\ref{sequence}.)
Then $\2M$ is a smooth
submanifold through $p$ and, if $\eps_k$ have been chosen
sufficiently rapidly converging to $0$, $\2M$ is tangent to $M$ of
infinite order at $p$ in view of \eqref{epsk} and \eqref{epsk1}. Consequently $\2M$
is also of finite type at $p$. Furthermore, the germ $(\2M,p)$ is clearly
invariant under the action of $G$, i.e.
the group $\Aut_\CR(\2M,p)$ of germs at $p$ of all local CR-automorphisms
of $M$ fixing $p$ contains $G$.

We now claim that $\Aut_\CR(\2M,p)=G$. Indeed, fix any $f\in
\Aut_\CR(\2M,p)$ and its representative defined in some neighborhood of $p$
in $\2M$, denoted by the same letter. Then for $k$ sufficiently
large, $f$ is defined in $D_k\cap \2M$ with $f(D_k\cap \2M)\subset D_1$.
By the construction, each $q_k\in N_k\cap U_k$ is not contained in $W_k$,
hence it is in $M_k\cap V_k$ and therefore in $\2M$
so that we can evaluate $f(q_k)$.
Then \eqref{orb-cond} implies that $f(q_k)\notin M$ and hence
$f(q_k)\in g(U_s)\subset g(V_s)$ for some $s$ and some $g\in G_s$.
Thus we have the estimates
\begin{equation}\Label{}
(\sup_{D_s} \|(g^{-1})'\|)^{-1}\, 
\frac{\inf_{z\in V_s} |z-p|}{\sup_{z\in V_k}|z-p|}
\le \frac{|f(q_k)-p|}{|q_k-p|} 
\le \sup_{D_s} \|g'\|\, 
\frac{\sup_{z\in V_s} |z-p|}{\inf_{z\in V_k}|z-p|}.
\end{equation}
On the other hand,
since $f$ is a local diffeomorphism of $\2M$ fixing $p$,
there exist constants $0<c<C$ such that 
\begin{equation}\Label{f-est}
c \le \frac{|f(z)-p|}{|z-p|} \le C, \quad c \le \|f'(z)\| \le C
\end{equation}
for $z\ne p\in \2M$ sufficiently close to $p$.
Then if $k$ is sufficiently large, we must have $s=k$ in view
of \eqref{supinf}. Hence $f(q_k)\in g_k(U_k)$ for suitable $g_k\in G_k$.
Setting $f_k:=g_k^{-1}\circ f\in \Aut_\CR(M,p)$, we have $f_k(q_k)\in U_k\cap \2M$. 
In view of \eqref{consm}, this means $f_k(q_k)\in U_k\cap M_k$. 

We now claim that, for $k$ sufficiently large,
we must have $B_{\eps_k\delta_k}(f_k(q_k)) \cap \2M\subset N_k$.
Indeed, otherwise, in view of \eqref{mk-nk},
we would have $W_k\cap B_{\eps_k\delta_k}(f_k(q_k))\ne \emptyset$
and, since $\eta_k < \eps_k\delta_k$, it would imply 
$\2X_k \cap B_{2\eps_k\delta_k}(f_k(q_k))\ne \emptyset$.
However, in view of \eqref{eta} and \eqref{f-est},
this would mean that, for $k$ sufficiently large
and some point $y\in \2X_k$, the inclusion 
$f_k^{-1}(y)\in B_{\delta_k}(q_k)\cap \2M$ would hold.
By our construction, $B_{\delta_k}(q_k)\cap \2M\subset N_k$
and hence we would have a contradiction with \eqref{orb-cond1}.
Hence we have $B_{\eps_k\delta_k}(f_k(q_k)) \cap \2M\subset N_k$ as claimed.
Thus $B_{\eps_k\delta_k}(f_k(q_k)) \cap \2M = B_{\eps_k\delta_k}(f_k(q_k)) \cap N_k$.
Again, using \eqref{eta} and \eqref{f-est},
we conclude that, for $k$ sufficiently large, 
$f_k^{-1}$ sends $B_{\eps_k\delta_k}(f_k(q_k)) \cap N_k$
into $B_{\delta_k}(q_k)\cap N_k$.
By our construction of the set $X_k\subset N_k$,
the latter conclusion means either $f_k(q_k)\in X_k$ or $f_k(q_k)=q_k$.
The first case is impossible in view of the first condition in \eqref{mk-nk}.
Hence we have $f_k(q_k)=q_k$
and, since a neighborhood of $q_k$ in $M_k$ is totally
rigid, this means $f_k\equiv\id$ in a neighborhood of $q_k$. 
Since $(\2M,p)$ is of
finite type, we have $f_k\equiv\id$ as germs at $p$,
and hence $f\equiv g_k\in G$ implying the desired conclusion.
\epf

\end{document}